\def\draft{n}
\documentclass[12pt]{amsart}
\usepackage[headings]{fullpage}
\usepackage{amssymb,epic,eepic,epsfig}
\usepackage{graphicx}
\usepackage{texdraw}
\usepackage{url}
\usepackage[bookmarks=true,%
    colorlinks=true,%
    linkcolor=blue,%
    citecolor=blue,%
    filecolor=blue,%
    menucolor=blue,%
    urlcolor=blue,%
    breaklinks=true]{hyperref}

\newtheorem{theorem}{Theorem}[section]
\theoremstyle{definition}

\newtheorem{lemma}[theorem]{Lemma}
\newtheorem{definition}[theorem]{Definition}

\newtheorem{conjecture}[theorem]{Conjecture}
\newtheorem{question}[theorem]{Question}

\def\printname#1{
        \if\draft y
                \smash{\makebox[0pt]{\hspace{-0.5in}
                        \raisebox{8pt}{\tt\tiny #1}}}
        \fi
}

\newcommand{\psdraw}[2]
         {\begin{array}{c} \hspace{-1.3mm}
        \raisebox{-4pt}{\epsfig{figure=draws/#1.eps,width=#2}}
        \hspace{-1.9mm}\end{array}}

\newlength{\standardunitlength}
\catcode`\@=11
\long\def\@makecaption#1#2{%
     \vskip 10pt

\setbox\@tempboxa\hbox{
       \small\sf{\bfcaptionfont #1. }\ignorespaces #2}%
     \ifdim \wd\@tempboxa >\captionwidth {%
         \rightskip=\@captionmargin\leftskip=\@captionmargin
         \unhbox\@tempboxa\par}%
       \else
         \hbox to\hsize{\hfil\box\@tempboxa\hfil}%
     \fi}
\font\bfcaptionfont=cmssbx10 scaled \magstephalf
\newdimen\@captionmargin\@captionmargin=2\parindent
\newdimen\captionwidth\captionwidth=\hsize
\catcode`\@=12

\def\lbl#1{\label{#1}\printname{#1}}

\def\BN{\mathbb N}
\def\BZ{\mathbb Z}

\def\BQ{\mathbb Q}
\def\BR{\mathbb R}
\def\BC{\mathbb C}

\def\calB{\mathcal B}

\def\D{\Delta}

\def\a{\alpha}

\def\l{\lambda}

\def\s{\sigma}
\def\ga{\gamma}

\def\la{\langle}
\def\ra{\rangle}

\def\e{\epsilon}

\def\d{\delta}
\def\b{\beta}

\def\s{\sigma}

\def\longto{\longrightarrow}

\def\SL{\mathrm{SL}}
\def\PSL{\mathrm{PSL}}
\def\GL{\mathrm{GL}}
\def\pt{\partial}

\def\calB{\mathcal{B}}
\def\z{\zeta}

\def\bs{\mathrm{bs}}
\def\fsl{\mathfrak{sl}}
\def\mindeg{\mathrm{mindeg}}

\def\be{\begin{equation}}
\def\ee{\end{equation}}

\def\Hom{\mathrm{Hom}}
\def\hb{\hbar}
\def\ve{\varepsilon}

\begin{document}


\title[Quantum Knot Invariants]{Quantum Knot Invariants}
\author{Stavros Garoufalidis}
\address{School of Mathematics \\
         Georgia Institute of Technology \\
         Atlanta, GA 30332-0160, USA \newline
         {\tt \url{http://www.math.gatech.edu/~stavros }}}
\email{stavros@math.gatech.edu}
\thanks{The author was supported in part by NSF. \\
\newline
1991 {\em Mathematics Classification.} Primary 57N10. Secondary 57M25.
}

\date{January 15, 2012}

\dedicatory{To Don Zagier, with admiration}

\begin{abstract}
This is a survey talk on one of the best known quantum knot 
invariants, the colored Jones polynomial of a knot, and its relation to
the algebraic/geometric topology and hyperbolic geometry of the 
knot complement. We review several aspects of the colored
Jones polynomial, emphasizing modularity, stability and effective computations.
The talk was given in the Mathematische Arbeitstagung June 24-July 1, 2011.
\end{abstract}

\maketitle

\tableofcontents

\section{The Jones polynomial of a knot}
\lbl{sec.jones}

Quantum knot invariants are powerful numerical invariants defined 
by Quantum Field theory with deep connections to the geometry and topology
in dimension three \cite{Wi}. This is a survey talk on the various limits
the colored Jones polynomial \cite{Jo}, one of the best known quantum knot 
invariants. This is a 25 years old subject that contains theorems and
conjectures in disconnected areas of mathematics. We chose to present 
some old and recent conjectures on the subject, emphasizing two recent 
aspects of the colored Jones polynomial, Modularity and Stability
and their illustration by effective computations. Zagier's influence on 
this subject is profound, and several results in this talk are joint 
work with him. Of course, the author is responsible for any mistakes in the 
presentation. We thank Don Zagier for enlightening conversations, for
his hospitality and for his generous sharing of his ideas with us.

The Jones polynomial $J_L(q) \in \BZ[q^{\pm 1/2}]$ of an oriented link $L$ in 
3-space is uniquely determined by the linear relations \cite{Jo}
$$
q J_{\psdraw{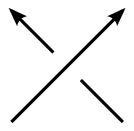}{0.3in}}(q) - q^{-1} J_{\psdraw{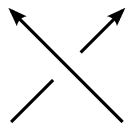}{0.3in}}(q)
=(q^{1/2}-q^{-1/2}) J_{\psdraw{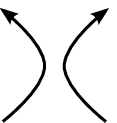}{0.3in}}(q)
\qquad
J_{\psdraw{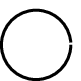}{0.3in}}(q)=q^{1/2}+q^{-1/2} \,.
$$
The Jones polynomial has a unique extension to a polynomial invariant
$J_{L,c}(q)$ of links $L$ together with a coloring $c$ of their 
components are colored by positive natural 
numbers that satisfy the following rules
\begin{eqnarray*}
J_{L \cup K, c \cup \{N+1\}}(q) &=& J_{L \cup K^{(2)}, c \cup \{N,2\} }(q)-
J_{L \cup K, c \cup \{N-1\}}(q), \qquad N \geq 2 \\
J_{L \cup K, c \cup \{1\}}(q) &=& J_{L, c}(q) \\
J_{L,\{2,\dots,2\}}(q) &=& J_L(q)
\end{eqnarray*}
where $(L \cup K, c \cup \{N\})$ denotes a link with a 
distinguished component $K$ colored by $N$
and $K^{(2)}$ denotes the 2-parallel of $K$ with zero framing.
Here, a natural number $N$ attached to a component of a link indicates the
$N$-dimensional irreducible representation of the Lie algebra $\fsl(2,\BC)$.
For a detailed discussion on the polynomial invariants of links that come
from quantum groups, see \cite{Ja,Tu1,Tu2}.

The above relations make clear that the colored Jones polynomial of a knot
encodes the Jones polynomials of the knot and its 0-framed parallels.



\section{Three limits of the colored Jones polynomial}
\lbl{sec.3limits}

In this section we will list three conjectures, the MMR Conjecture (proven),
the Slope Conjecture (mostly proven) and the AJ Conjecture (less proven).
These conjectures relate the colored Jones polynomial
of a knot with the Alexander polynomial, with the set of slopes of 
incompressible surfaces and with the $\PSL(2,\BC)$ character variety
of the knot complement.

\subsection{The colored Jones polynomial and the Alexander polynomial}
\lbl{sub.MMR}

We begin by discussing a relation of the colored Jones polynomial
of a knot with the homology of the universal abelian cover of its complement.
The homology $H_1(M,\BZ)\simeq\BZ$ of the complement $M=S^3\setminus K$ of a 
knot $K$ in 3-space is independent of the knot $K$. This allows us to
consider the universal abelian cover $\widetilde{M}$ of $M$ with deck
transformation group $\BZ$, and with homology $H_1(\widetilde{M},\BZ)$
a $\BZ[t^{\pm 1}]$ module. As is well-known this module is essentially 
torsion and its order is given by the Alexander polynomial 
$\D_K(t) \in \BZ[t^{\pm 1}]$ of $K$ \cite{Rf}. The Alexander polynomial does
not distinguish knots from their mirrors and satisfies $\D_K(1)=1$.

There are infinitely many pairs of knots (for instance $(10_{22},10_{35})$
in the Rolfsen table \cite{Rf,B-N}) with equal 
Jones polynomial but different Alexander polynomial. On the other hand, 
the colored Jones polynomial 
determines the Alexander polynomial. This so-called Melvin-Morton-Rozansky
Conjecture was proven in \cite{B-NG}, and states that
\begin{equation}
\lbl{eq.MMR}
\hat J_{K,n}(e^{\hb})
=\sum_{i \geq j \geq 0} a_{K,ij} \hb^i n^j \in \BQ[[n,\hb]]
\end{equation}
and
$$
\sum_{i=0}^\infty a_{K,ij} \hb^i=\frac{1}{\D_K(e^{\hb})} \in \BQ[[\hb]]\,.
$$
Here $\hat J_{K,n}(q)=J_{K,n}(q)/J_{\text{Unknot},n}(q) \in \BZ[q^{\pm 1}]$ 
is a normalized form of the colored Jones polynomial. The above conjecture
is a statement about formal power series. A stronger analytic version
is known \cite[Thm.1.3]{GL2}, namely for every knot $K$ there exists
an open neighborhood $U_K$ of $0 \in \BC$ such that for all $\a \in U_K$
we have
$$
\lim_n J_{K,n}(e^{\a/n})=\frac{1}{\D_K(e^{\a})}\,,
$$
where convergence is uniform with respect to compact sets. More is known
about the summation of the series \eqref{eq.MMR} along a fixed diagonal
$i=j+k$ for fixed $k$, both on the level of formal power series and on
the analytic counterpart. For further details the reader may consult
\cite{GL2} and references therein.

\subsection{The colored Jones polynomial and slopes of incompressible
surfaces}
\lbl{sub.slope}

In this section we discuss a conjecture relating the degree
of the colored Jones polynomial of a knot $K$ with the set $\bs_K$ 
of boundary slopes of incompressible surfaces in the knot complement 
$M=S^3\setminus K$. Although there are infinitely many incompressible
surfaces in $M$, it is known that $\bs_K \subset \BQ \cup \{1/0\}$ is
a finite set \cite{Ht}. Incompressible surfaces play an important role
in geometric topology in dimension three, often accompanied by the theory
of normal surfaces \cite{Hk}. 
From our point of view, incompressible surfaces are
a tropical limit of the colored Jones polynomial, corresponding to an
expansion around $q=0$ \cite{Ga5}.

The Jones polynomial of a knot is a Laurent polynomial in one variable $q$
with integer coefficients. Ignoring most information, one can consider the
degree $\d_K(n)$ of $\hat J_{K,n+1}(q)$ with respect to $q$. 
Since $(\hat J_{K,n}(q))$ is a $q$-holonomic sequence
\cite{GL}, it follows that
$\d_K$ is a quadratic quasi-polynomial \cite{Ga4}. In other words, we have
$$
\d_K(n)=c_K(n) n^2 + b_K(n) n + a_K(n)\,,
$$
where $a_K, b_K, c_K: \BN \longto \BQ$ are periodic functions. In \cite{Ga3}
the author formulated the Slope Conjecture.

\begin{conjecture}
\lbl{conj.slope}
For all knots $K$ we have
$$
4 c_K(\BN) \subset \bs_K \,.
$$
\end{conjecture}
The movitating example for the Slope Conjecture was
the case of the $(-2,3,7)$ pretzel knot, where we have
\cite[Ex.1.4]{Ga3}
$$
\d_{(-2,3,7)}(n) = \left[\frac{37}{8} n^2 + \frac{17}{2} n \right]=
\frac{37}{8} n^2 + \frac{17}{2} n + a(n),
$$
where $a(n)$ is a periodic sequence of period $4$ given by $0,-1/8,-1/2,-1/8$
if $n\equiv 0,1,2,3 \bmod 4$ respectively. In addition, we have
$$
\bs_{(-2,3,7)}=\{0,16, 37/2, 20\} \,.
$$
In all known examples, $c_K(\BN)$ consists of a single element, the
so-called Jones slope. How the colored Jones polynomial selects some of the
finitely many boundary slopes is a challenging and interesting question.
The Slope Conjecture is known for all torus knots, all alternating knots and 
all knots with at most $8$ crossings \cite{Ga3} as well as for all 
adequate knots \cite{FKP} and all 2-fusion knots \cite{DunG}.

\subsection{The colored Jones polynomial and the $\PSL(2,\BC)$
character variety}
\lbl{sub.AJ}

In this section we discuss a conjecture relating the colored 
Jones polynomial of a knot $K$ with the moduli space of 
$\SL(2,\BC)$-representations of $M$, restricted to the boundary of $M$.
Ignoring $0$-dimensional components, the latter is a 1-dimensional 
plane curve. To formulate the conjecture we need to recall that the
colored Jones polynomial $\hat J_{K,n}(q)$ is $q$-holonomic \cite{GL} i.e.,
it satisfies a non-trivial linear recursion relation
\begin{equation}
\lbl{eq.recJ}
\sum_{j=0}^d a_j(q,q^n) \hat J_{K,n+j}(q)=0
\end{equation}
for all $n$ where $a_j(u,v) \in \BZ[u^{\pm 1},v^{\pm 1}]$ and $a_d \neq 0$.
$q$-holonomic sequences were introduced by Zeilberger \cite{Zeil}, and
a fundamental theorem (multisums of $q$-proper hypergeometric terms are
$q$-holonomic) was proven in \cite{WZ} and implemented in \cite{PWZ}.
Using two operators $M$ and $L$ which act on a sequence $f(n)$ by
$$
(M f)(n)=q^n f(n), \qquad (L f)(n)=f(n+1) \,,
$$
we can write the recursion \eqref{eq.recJ} in operator form 
$$
P \cdot \hat J_{K}=0 \qquad \text{where} \qquad P=\sum_{j=0}^d a_j(q,M) L^j \,.
$$
It is easy to see that $LM=qML$ and $M,L$ generate the $q$-Weyl algebra.
One can choose a canonical 
recursion $A_K(M,L,q) \in \BZ[q,M]\la L \ra/(LM-qML)$ which is a knot
invariant \cite{Ga2}, the non-commutative $A$-polynomial of $K$.
The reason for this terminology is the potential relation with the
$A$-polynomial $A_K(M,L)$ of $K$ \cite{CCGLS}. The latter is defined
as follows. 

Let $X_M=\Hom(\pi_1(M),\SL(2,\BC))/\BC$ denote the moduli space of
flat $\SL(2,\BC)$ connections on $M$. We have an identification
$$
X_{\pt M} \simeq (\BC^*)^2/(\BZ/2\BZ), \qquad \rho \mapsto (M,L)
$$
where $\{M,1/M\}$ (resp., $\{L,1/L\}$) are the eigenvalues of $\rho(\mu)$
(resp., $\rho(\l)$) where $(\mu,\l)$ is a meridian-longitude pair on
$\pt M$.
$X_M$ and $X_{\pt M}$ are affine varieties and the restriction map
$X_M \longto X_{\pt M}$ is algebraic. The Zariski closure of its image
lifted to $(\BC^*)^2$, and after removing any $0$-dimensional components
is a one-dimensional plane curve with defining polynomial $A_K(M,L)$
\cite{CCGLS}. This polynomial plays an important role in the hyperbolic 
geometry of the knot complement. We are now ready to formulate the
AJ Conjecture \cite{Ga2}; see also \cite{Ge}. Let us say that two 
polynomials $P(M,L)=_M Q(M,L)$ are essentially equal if their irreducible 
factors with positive $L$-degree are equal.

\begin{conjecture}
\lbl{conj.AJ}
For all knots $K$, we have $A_K(M^2,L,1)=_M A_K(M,L)$.
\end{conjecture}
The AJ Conjecture was checked for the $3_1$ and the $4_1$ knots in \cite{Ga2}.
It is known for most 2-bridge knots \cite{Le}, for torus knots and for the 
pretzel knots of Section \ref{sec.computeAhat}; see \cite{LT,Tran}.

From the point of view of physics, the AJ Conjecture is a consequence
of the fact that quantization and the corresponding quantum field theory
exists \cite{Gu,Di}.

\section{The Volume and Modularity Conjectures}
\lbl{sec.VC}

\subsection{The Volume Conjecture}
\lbl{sub.VC}

The Kashaev invariant of a knot is a sequence of complex numbers defined by
\cite{Ka,MM}
$$
\la K \ra_N=\hat J_{K,N}(e(1/N))
$$
where $e(\a)=e^{2 \pi i \a}$. The Volume Conjecture concerns the exponential
growth rate of the Kashaev invariant and states that
$$
\lim_N \frac{1}{N} \log|\la K \ra_N|=\frac{\text{vol}(K)}{2\pi}
$$
where $\text{Vol(K)}$ is the volume of the hyperbolic pieces of the knot
complement $S^3\setminus K$ \cite{Th}. 
Among hyperbolic knots, the Volume Conjecture
is known only for the $4_1$ knot. Detailed computations are available
in \cite{Mu1}. 
Refinements of the Volume Conjecture to all orders in $N$ and generalizations 
were proposed by several authors \cite{DGLZ,GuMu,GL2,Ga0}. Although proofs are 
lacking, there appears to be a lot of structure in the asymptotics of the
Kashaev invariant. In the next section we will discuss a modularity 
conjecture of Zagier and some numerical verification.

\subsection{The Modularity Conjecture}
\lbl{sec.modularity}

Zagier considered the Galois invariant spreading of the Kashaev invariant
on the set of complex roots of unity given by
$$
\phi_K: \BQ/\BZ \longto \BC,
\qquad
\phi_K\left(\frac{a}{c}\right)=\hat J_{K,c}\left(e\left(\frac{a}{c}\right)\right)
$$
where $(a,c)=1$ and $c>0$. The above formula works even when $a$ and $c$
are not coprime due to a symmetry of the colored Jones
polynomial \cite{Ha}. $\phi_K$ determines $\la K \ra$ and 
conversely is determined by $\la K \ra$ via Galois invariance. 

Let $\ga=\left(\begin{matrix} a & b \\ c & d \end{matrix}\right) \in 
\SL(2,\BZ)$ and $\a=a/c$ and $\hb=2 \pi i/(X+d/c)$ where $X \longto +\infty$ 
with bounded denominators. Let $\phi=\phi_K$ denote the extended Kashaev
invariant of a hyperbolic knot $K$ and let $F \subset \BC$ denote the 
invariant trace field of $M=S^3\setminus K$ \cite{MR}. 
Let $C(M) \in \BC/(4 \pi^2 \BZ)$ denote 
the complex Chern-Simons invariant of $M$ \cite{GoZi,Ne}. The next
conjecture was formulated by Zagier.

\begin{conjecture}
\lbl{conj.modular}\cite{Za:QMF}
With the above conventions, there exist $\D(\a) \in \BC$ with 
$\D(\a)^{2c} \in F(\e(\a))$ and $A_j(\a) \in F(e(\a))$ such that
\begin{equation}
\lbl{eq.Kgamma}
\frac{\phi(\ga X)}{\phi(X)} \sim \left(\frac{2\pi}{\hb}\right)^{3/2}
e^{C(M)/\hb} \D(\a) \sum_{j=0}^\infty A_j(\a) \hb^j \,.
\end{equation}
\end{conjecture}
When $\ga=\left(\begin{matrix} 1 & 0 \\ 1 & 1 \end{matrix}\right)$
and $X=N-1$, and with the properly chosen orientation of $M$,
the leading asymptotics of \eqref{eq.Kgamma} together
with the fact that $\Im(C(M))=\text{vol}(M)$ gives the volume conjecture.

\section{Computation of the non-commutative $A$-polynomial}
\lbl{sec.computeAhat}

As we will discuss below, the key to an effective computation 
the Kashaev invariant is a recursion for the colored Jones polynomial.
Proving or guessing such a recursion is at least as hard as computing 
the $A$-polynomial of the knot.
The $A$-polynomial is already unknown for several knots with $9$ crossings.
For an updated table of computations see \cite{Cu}. The $A$-polynomial is
known for the 1-parameter
families of twist knots $K_p$ \cite{HS} and pretzel knots 
$KP_p=(-2,3,3+2p)$ \cite{GM} 
depicted on the left and the right part of the following figure
$$
\psdraw{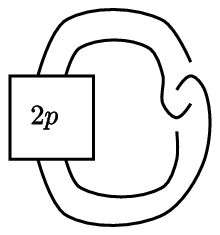}{0.9in} \qquad\qquad \psdraw{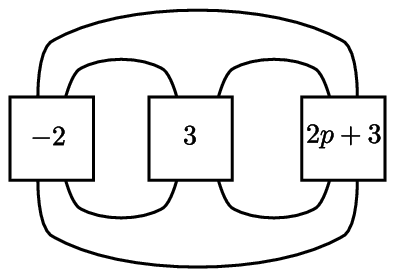}{1.5in}
$$
where an integer $m$ inside a box indicates the number 
of $|m|$ half-twists, right-handed (if $m>0$) or left-handed (if $m<0$), 
according to the following figure
$$
\psdraw{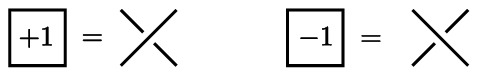}{2in}
$$
The non-commutative $A$-polynomial of the twist knots $K_p$ was computed 
with a certificate by X. Sun and the author in \cite{GS2} 
for $p=-14,\dots,15$. The data is available from
\begin{center}
\url{http://www.math.gatech.edu/~stavros/publications/twist.knot.data}
\end{center}
The non-commutative $A$-polynomial of the pretzel knots $KP_p=(-2,3,3+3p)$
was guessed by C. Koutschan and the author in  
\cite{GK} for $p=-5,\dots,5$. The guessing method used an a priori
knowledge of the monomials of the recursion, together with computation of
the colored Jones polynomial using the fusion formula, and exact but modular
arithmetic and rational reconstruction. The data is available  from
\begin{center}
\url{http://www.math.gatech.edu/~stavros/publications/pretzel.data}
\end{center}
For instance, the recursion relation for the colored Jones polynomial
$f(n)$ of the $5_2=(-2,3,-1)$ pretzel knot is given by 
\vspace{0.2cm}


{\small
\begin{math}
b(q^n,q)
-q^{9+7 n} (-1+q^n) (-1+q^{2+n}) (1+q^{2+n}) (-1+q^{5+2 n}) f(n)+q^{5+2 n} (-1+q^{1+n})^2 (1+q^{1+n}) (-1+q^{5+2 n}) (-1+q^{1+n}+q^{1+2 n}-q^{2+2 n}-q^{3+2 n}+q^{4+2 n}-q^{2+3 n}-q^{5+3 n}-2 q^{5+4 n}+q^{6+5 n}) f(1+n)-q (-1+q^{2+n})^2 (1+q^{2+n}) (-1+q^{1+2 n}) (-1+2 q^{2+n}+q^{2+2 n}+q^{5+2 n}-q^{4+3 n}+q^{5+3 n}+q^{6+3 n}-q^{7+3 n}-q^{7+4 n}+q^{9+5 n}) f(2+n)-(-1+q^{1+n}) (1+q^{1+n}) (-1+q^{3+n}) (-1+q^{1+2 n}) f(3+n)=0 \,,
\end{math}
}

\vspace{0.2cm}
where
\vspace{0.2cm}

{\small
\begin{math}
b(q^n,q)=q^{4+2 n} (1+q^{1+n}) (1+q^{2+n}) (-1+q^{1+2 n}) (-1+q^{3+2 n}) (-1+q^{5+2 n}) \,.
\end{math}
}

\vspace{0.2cm}
The recursion relation for the colored Jones polynomial $f(n)$ of the 
$(-2,3,7)$ pretzel knot is given by
\vspace{0.2cm}

{\small
\begin{math}
b(q^n,q)
-q^{224+55 n} (-1+q^n) (-1+q^{4+n}) (-1+q^{5+n}) f(n)+q^{218+45 n} (-1+q^{1+n})^3 (-1+q^{4+n}) (-1+q^{5+n}) f(1+n)+q^{204+36 n} (-1+q^{2+n})^2 (1+q^{2+n}+q^{3+n}) (-1+q^{5+n}) f(2+n)+(-1+q) q^{180+27 n} (1+q) (-1+q^{1+n}) (-1+q^{3+n})^2 (-1+q^{5+n}) f(3+n)-q^{149+18 n} (-1+q^{1+n}) (-1+q^{4+n})^2 (1+q+q^{4+n}) f(4+n)-q^{104+8 n} (-1+q^{1+n}) (-1+q^{2+n}) (-1+q^{5+n})^3 f(5+n)+q^{59} (-1+q^{1+n}) (-1+q^{2+n}) (-1+q^{6+n}) f(6+n)=0\,,
\end{math}
}

\vspace{0.2cm}
where
\vspace{0.2cm}

{\small
\begin{math}
b(q^n,q)=
q^{84+5 n} (1-q^{1+n}-q^{2+n}+q^{3+2 n}-q^{16+3 n}+q^{17+4 n}+q^{18+4 n}-q^{19+5 n}-q^{26+5 n}+q^{27+6 n}+q^{28+6 n}+q^{31+6 n}-q^{29+7 n}-q^{32+7 n}-q^{33+7 n}-q^{36+7 n}+q^{34+8 n}+q^{37+8 n}+q^{38+8 n}-q^{39+9 n}+q^{45+9 n}-q^{46+10 n}-q^{47+10 n}+q^{49+10 n}+q^{48+11 n}-q^{50+11 n}-q^{51+11 n}-q^{54+11 n}+q^{52+12 n}+q^{55+12 n}+q^{56+12 n}-q^{57+13 n}-q^{62+13 n}+q^{63+14 n}+q^{64+14 n}-q^{66+14 n}+q^{67+14 n}-q^{65+15 n}+q^{67+15 n}-q^{69+15 n}+q^{71+15 n}-q^{69+16 n}+q^{70+16 n}-q^{72+16 n}-q^{75+17 n}-q^{78+17 n}+q^{76+18 n}+q^{79+18 n}-q^{83+19 n}+q^{85+19 n}+q^{84+20 n}-q^{86+20 n}+q^{88+20 n}-q^{89+21 n}+q^{91+21 n}-q^{96+22 n}-q^{93+23 n}+2 q^{98+24 n}-q^{99+25 n}-q^{108+26 n}-q^{107+27 n}+q^{109+27 n}+q^{108+28 n}-q^{110+28 n}+q^{112+28 n}-q^{113+29 n}+q^{115+29 n}+q^{112+30 n}+q^{115+30 n}-q^{117+31 n}-q^{120+31 n}-q^{117+32 n}+q^{118+32 n}-q^{120+32 n}-q^{119+33 n}+q^{121+33 n}-q^{123+33 n}+q^{125+33 n}+q^{123+34 n}+q^{124+34 n}-q^{126+34 n}+q^{127+34 n}-q^{123+35 n}-q^{128+35 n}+q^{124+36 n}+q^{127+36 n}+q^{128+36 n}+q^{126+37 n}-q^{128+37 n}-q^{129+37 n}-q^{132+37 n}-q^{130+38 n}-q^{131+38 n}+q^{133+38 n}-q^{129+39 n}+q^{135+39 n}+q^{130+40 n}+q^{133+40 n}+q^{134+40 n}-q^{131+41 n}-q^{134+41 n}-q^{135+41 n}-q^{138+41 n}+q^{135+42 n}+q^{136+42 n}+q^{139+42 n}-q^{133+43 n}-q^{140+43 n}+q^{137+44 n}+q^{138+44 n}-q^{142+45 n}+q^{135+46 n}-q^{139+47 n}-q^{140+47 n}+q^{144+48 n})\,.
\end{math}
}

\vspace{0.2cm}
The pretzel knots $KP_p$ are interesting from many points of view.
For every integer $p$, the knots in the pair $(KP_p, -KP_{-p})$ 
(where $-K$ denotes the mirror of $K$)
\begin{itemize}
\item
are geometrically similar, in particular they are scissors congruent,
have equal volume, equal invariant trace fields and   
their Chern-Simons invariant differ by a sixth root of unity,
\item
their $A$-polynomials are equal up to a $\GL(2,\BZ)$ transformation
\cite[Thm.1.4]{GM}.
\end{itemize}

Yet, the colored Jones polynomials of $(KP_p, -KP_{-p})$ are different, 
and so are the Kashaev invariants and their asymptotics and even the 
term $\D(0)$ in the modularity conjecture \ref{conj.modular}.
An explanation of this difference is given in \cite{DG}.

Zagier posed a question to compare the modularity conjecture for geometrically
similar pairs of knots, which was a motivation for many of the computations
in Section \ref{sub.mnum}.

\section{Numerical asymptotics and the Modularity Conjecture}
\lbl{sec.num}

\subsection{Numerical computation of the Kashaev invariant}
\lbl{sub.knum}

To numerically verify Conjecture \ref{conj.modular} we need to compute
the Kashaev invariant to several hundred digits when $N=2000$
for instance. In this section we discuss how to achieve this.

There are multidimensional $R$-matrix 
state sum formulas for the colored Jones polynomial $J_{K,N}(q)$ where the
number of summation points are given by a polynomial in $N$ of degree
the number of crossings of $K$ minus $1$ \cite{GL}. Unfortunately, this
is not practical method even for the $4_1$ knot. 

An alternative way is to use fusion \cite{KL,Co,GV} which allows one to 
compute the colored Jones polynomial more efficiently at the cost that the 
summand is a rational function of $q$. For instance, the colored Jones 
polynomial of a 2-fusion knot can be computed in 
$O(N^3)$ steps using \cite[Thm.1.1]{GK}. This method works better, but
it still has limitations.

A preferred method is to guess a nontrivial recursion relation for the 
colored Jones polynomial (see Section \ref{sec.computeAhat}) and instead
of using it to compute the colored Jones polynomial, differentiate 
sufficiently many times and numerically compute the Kashaev invariant.
In the efforts to compute the Kashaev invariant of the $(-2,3,7)$ pretzel
knot, Zagier and the author obtained the following lemma, of
theoretical and practical use.

\begin{lemma}
\lbl{lem.1}
The Kashaev invariant $\la K \ra_N$ can be numerically computed in $O(N)$
steps.
\end{lemma}
A computer implementation of Lemma \ref{lem.1} is available.

\subsection{Numerical verification of the  Modularity Conjecture}
\lbl{sub.mnum}

Given a sequence of complex number $(a_n)$ with an expected asymptotic
expansion
$$
a_n \sim \l^n n^{\a} (\log n)^{\b} \sum_{j=0}^\infty \frac{c_j}{n^j}
$$ 
how can one numerically compute $\l$, $\a$, $\b$ and several coefficients
$c_j$? This is a well-known numerical analysis problem \cite{BO}. 
An acceleration
method was proposed in \cite[p.954]{Za2}, which is also equivalent to the
Richardson transform. For a detailed discussion of the acceleration
method see \cite[Sec.5.2]{GIKM}. In favorable circumstances the coefficients
$c_j$ are algebraic numbers, and a numerical approximation may lead to a guess
for their exact value.

A concrete application of the acceleration method was given in the appendix 
of \cite{GV} where one deals with several $\l$ of the same magnitude as well 
as $\b \neq 0$. 

Numerical computations of the modularity conjecture for the $4_1$
knot were obtained by Zagier around roots of unity of order at most $5$,
and extended to several other knots in \cite{GZ1,GZ2}. As a sample computation,
we present here the numerical data for $4_1$ at $\a=0$, 
computed independently by Zagier and by the author. The values of
$A_k$ in the table below are known for $k=0,\dots,150$.

$$
\phi_{4_1}(X) = 3^{-1/4} X^{3/2} \exp(C X) \left(
\sum_{k=0}^\infty \frac{A_k}{k! 12^k} h^k \right)
$$
$$
h=A/X \qquad A=\frac{\pi}{3^{3/2}} 
\qquad C=\frac{1}{\pi}\mathrm{Li_2}(\exp(2 \pi i/3))
$$

$$
\begin{array}{|r|l|} \hline
k & A_k \\ \hline
0 & 1 \\ \hline 
1 & 11 \\ \hline 
2 & 697 \\ \hline 
3 & 724351/5 \\ \hline 
4 & 278392949/5 \\ \hline 
5 & 244284791741/7 \\ \hline 
6 & 1140363907117019/35 \\ \hline 
7 & 212114205337147471/5 \\ \hline
8 & 367362844229968131557/5 \\ \hline 
9 & 44921192873529779078383921/275 \\ \hline 
10 & 3174342130562495575602143407/7 \\ \hline
11 & 699550295824437662808791404905733/455 \\ \hline 
12 & 14222388631469863165732695954913158931/2275 \\ \hline
13 & 5255000379400316520126835457783180207189/175 \\ \hline 
14 & 4205484148170089347679282114854031908714273/25 \\ \hline
15 & 16169753990012178960071991589211345955648397560689/14875 \\ \hline 
16 & 119390469635156067915857712883546381438702433035719259/14875 \\ \hline
17 & 1116398659629170045249141261665722279335124967712466031771/16625 \\ \hline
18 & 577848332864910742917664402961320978851712483384455237961760783/914375 \\ \hline
19 & 319846552748355875800709448040314158316389207908663599738774271783/48125 \\ \hline
20 & 5231928906534808949592180493209223573953671704750823173928629644538303/67375 \\ \hline
\end{array}
$$
{\tiny
$$
\begin{array}{|r|l|} \hline
21 & 158555526852538710030232989409745755243229196117995383665148878914255633279/158125 \\ \hline
22 & 2661386877137722419622654464284260776124118194290229321508112749932818157692851/186875 \\ \hline
23 & 1799843320784069980857785293171845353938670480452547724408088829842398128243496119/8125 \\ \hline
24 & 1068857072910520399648906526268097479733304116402314182132962280539663178994210946666679/284375 \\ \hline
25 & 1103859241471179233756315144007256315921064756325974253608584232519059319891369656495819559/15925 \\ \hline
26 & 8481802219136492772128331064329634493104334830427943234564484404174312930211309557188151604709/6125 \\ \hline
\end{array}
$$
}

In addition, we present the numerical data for the $5_2$ knot 
at $\a=1/3$, computed in \cite{GZ1}. 


\begin{center}
\end{center}
$$
\phi_{5_2}(X/(3X+1))/\phi_{5_2}(X) \sim e^{C/h} (2 \pi/h)^{3/2} \D(1/3) 
\left( \sum_{k=0}^\infty A_k(1/3) h^k  \right)
$$
$$h=(2 \pi i)/(X+1/3)$$ 
$$ 
F=\BQ(\a) \qquad \a^3-\a^2+1=0 \qquad \a=0.877\dots-0.744\dots i
$$
$$
C=R(1-\a^2) + 2 R(1-\a) -\pi i \log(\a)+\pi^2
$$
$$
R(x)=\mathrm{Li}_2(x)+\frac{1}{2} \log x \log(1-x)-\frac{\pi^2}{6}
$$
$$
[1-\a^2]+2[1-\a] \in \calB(F)
$$
$$
-23 = \pi_1^2 \pi_2 \qquad \pi_1=3\a-2 \qquad \pi_2=3\a+1
$$
$$
\pi_7=(\a^2-1)\z_6-\a+1 \qquad \pi_{43}=2\a^2-\a-\z_6
$$
\begin{eqnarray*}
\D(1/3) &=& e(-2/9) \pi_7 \frac{3 \sqrt{-3}}{\sqrt{\pi_1}} \\
A_0(1/3) &=& \pi_7 \pi_{43} \\
A_1(1/3) &=& \frac{-952 + 321 \a-873 \a^2 + (1348 + 557 \a + 26 \a^2)\z_6}{
\a^5 \pi_1^3}
\end{eqnarray*}

One may use the recursion relations \cite{GK2} for the twisted colored 
Jones polynomial to expand the above computations around complex roots 
of unity \cite{DG2}.

\section{Stability}
\lbl{sec.stability}

\subsection{Stability of a sequence of polynomials}
\lbl{sub.stables}

The Slope Conjecture deals with the highest (or the lowest, if you take
the mirror image) $q$-exponent of the
colored Jones polynomial. In this section we discuss what happens when
we shift the colored Jones polynomial and place its lowest $q$-exponent to 
$0$. Stability concerns the coefficients of the resulting sequence
of polynomials in $q$. A weaker form of stability ($0$-stability, defined 
below) for the colored Jones polynomial of an alternating knot was 
conjectured by Dasbach and Lin, and proven independently by Armond \cite{Ar2}.

Stability was observed in some examples of alternating knots by Zagier, and 
conjectured by the author to hold for all knots, assuming that we restrict 
the sequence of colored Jones polynomials to suitable arithmetic progressions, 
dictated by the quasi-polynomial nature of its $q$-degree \cite{Ga3,Ga4}. 
Zagier asked about modular and asymptotic properties of the limiting 
$q$-series.

A proof of stability in full for all alternating links is given in \cite{GL3}. 
Besides stability, this approach gives a generalized Nahm sum formula for the 
corresponding series, which in particular implies convergence in the open 
unit disk in the $q$-plane. The generalized Nahm sum
formula comes with a computer implementation (using as input a planar
diagram of a link), and allows the computation of several terms of the
$q$-series as well as its asymptotics when $q$ approaches
radially from within the unit circle a complex root of unity. 
The Nahm sum formula is reminiscent to the
cohomological Hall algebra of motivic Donaldson-Thomas invariants
of Kontsevich-Soibelman \cite{KS}, and may be related to recent work
of Witten \cite{Wi2} and Dimofte-Gaiotto-Gukov \cite{DGG}.

Let 
$$
\BZ((q ))=\{\sum_{n \in \BZ} a_n q^n \, | \, a_n=0, \, n \ll 0 \}
$$
denote the ring of power series in $q$ with integer coefficients and
bounded below minimum $q$-degree.

\begin{definition}
\lbl{def.stable}
Fix a sequence $(f_n(q))$ of polynomials $f_n(q) \in \BZ[q]$.
We say that  $(f_n(q))$ is $0$-{\em stable} if the following limit exists
\begin{equation*}
\lim_n f_n(q)=\Phi_0(q) \in \BZ[[q ]],
\end{equation*}
i.e.  for every natural
number $m\in \BZ$, there exists a natural number $n(m)$ such that the 
coefficient of $q^m$ in $f_n(q)$ is constant for all $n>n(m)$.

We say that $(f_n(q))$ is {\em stable} if there exist elements
$\Phi_k(q) \in \BZ((q))$ for $k=0,1,2,\dots$ such that for every
$k \in \BN$  we have
\begin{equation*}
\lim_n \, q^{-n k}\left(f_n(q)-\sum_{j=0}^{k} q^{jn} \Phi_j(q)\right)=0
\in \BZ((q))\,.
\end{equation*}
We will denote by
\begin{equation*}
F(x,q)=\sum_{k=0}^\infty \Phi_k(q) x^k \in \BZ((q ))[[x]]
\end{equation*}
the corresponding series associated to the stable sequence $(f_n(q))$.
%
\end{definition}

Thus, a $0$-stable sequence $f_n(q) \in \BZ[q]$ gives rise to a $q$-series
$\lim_n f_n(q) \in \BZ[[q]]$.  The $q$-series that come from the colored Jones
polynomial are $q$-hypergeometric series of a special shape, i.e., they are 
generalized Nahm sums. The latter are introduced in the next section.

\subsection{Generalized Nahm sums}
\lbl{sub.nahm}

In \cite{Nahm0} Nahm studied $q$-hypergeometric series $f(q) \in \BZ[[q]]$
of the form
\begin{equation*}
f(q) =\sum_{n_1,\dots,n_r \geq 0} \frac{q^{\frac{1}{2} n^t \cdot A \cdot n + b \cdot n}}{
(q)_{n_1} \dots (q)_{n_r}}
\end{equation*}
where $A$ is a positive definite even integral symmetric matrix and 
$b \in \BZ^r$.
Nahm sums appear in character formulas in Conformal Field Theory, and
define analytic functions in the complex unit disk $|q|<1$ with
interesting asymptotics at complex roots of unity, and with sometimes
modular behavior.  Examples of Nahm sums is the famous list of seven
mysterious $q$-series of Ramanujan that are nearly modular (in modern terms,
mock modular). For a detailed discussion, see \cite{Za.mock}. Nahm sums give
rise to elements of the Bloch group, which governs the leading radial
asymptotics of $f(q)$ as $q$ approaches a complex root of unity.
Nahm's Conjecture concerns the modularity of a Nahm sum $f(q)$, and was
studied extensively by Zagier, Vlasenko-Zwegers and others
\cite{VZ,Za.dilog}.

The limit of the colored Jones function of an alternating link
leads us to consider generalized Nahm sums of the form
\be
\lbl{eq.nahmgen}
\Phi(q)=\sum_{ n \in C \cap \BN^r} (-1)^{c \cdot n}
\frac{q^{\frac{1}{2} n^t \cdot A \cdot n + b \cdot n}}{(q)_{n_1} 
\dots (q)_{n_r}}
\ee
where $C$ is a rational polyhedral cone in $\BR^r$, $b, c \in \BZ^r$
and $A$ is a symmetric (possibly indefinite) symmetric matrix.
We will say that the generalized Nahm sum \eqref{eq.nahmgen} is 
{\em regular} if the function
$$
n \in C \cap \BN^r \mapsto \frac{1}{2} n^t \cdot A \cdot n + b \cdot n
$$
is proper and bounded below, where $\mindeg_q$ denotes the minimum degree
with respect to $q$. Regularity ensures that the series \eqref{eq.nahmgen}
is a well-defined element of the ring $\BZ((q))$.
In the remaining of the paper, the term Nahm 
sum will refer to a regular generalized Nahm sum.

\subsection{Stability for alternating links}
\lbl{sub.alts}

Let $K$ denote an alternating link. The lowest monomial of $J_{K,n}(q)$
has coefficient $\pm 1$, and dividing $J_{K,n+1}(q)$ by its lowest monomial
gives a polynomial $J^+_{K,n}(q) \in 1 + q \BZ[q]$.
We can now quote the main theorem of \cite{GL3}.

\begin{theorem}
\lbl{thm.2}\cite{GL3}
For every alternating link $K$, the sequence $(J^+_{K,n}(q))$ 
is stable and the corresponding limit $F_K(x,q)$ can be effectively computed
by a planar projection $D$ of $K$. Moreover, $F_K(0,q)=\Phi_{K,0}(q)$
is given by an explicit Nahm sum computed by $D$.
\end{theorem}
An illustration of the corresponding $q$-series $\Phi_{K,0}(q)$ the knots
$3_1$, $4_1$ and $6_3$ is given in Section \ref{sub.alt}.

\subsection{Computation of the $q$-series of alternating links}
\lbl{sub.alt}

Given the generalized Nahm sum for $\Phi_{K,0}(q)$, a multidimensional
sum of as many variables as the number of crossings of $K$, one may
try to identify the $q$-series $\Phi_{K,0}(q)$ with a known one. 
In joint work with Zagier,
we computed the first few terms of the corresponding series (an
interesting and nontrivial task in itself) and guessed the answer for knots 
with a small number of crossings. The guesses are presented in the following
table

{\small
$$
\begin{array}{|c|c|c|c|l|l|l|l|} \hline
K & c_- & c_+ & \s & \Phi^*_{K,0}(q) & \Phi_{K,0}(q) \\ \hline
\hline
3_1=-K_1 & 3 & 0 & 2 & h_3 & h_2 \\ \hline \hline
4_1=K_{-1} & 2 & 2 & 0 & h_3 & h_3 \\ \hline \hline
5_1 & 5 & 0 & 4 & h_5 & h_2 \\ \hline
5_2=K_2 & 0 & 5 & -2 & h_4 & h_3 \\ \hline \hline
6_1=K_{-2} & 4 & 2 & 0 & h_3 & h_5  \\ \hline
6_2 & 4 & 2 & 2 &  h_3 h_4 & h_3 \\ \hline
6_3 & 3 & 3 & 0 & h_3^2 & h_3^2 \\ \hline \hline
7_1 & 7 & 0 & 6 & h_7 & h_2 \\ \hline
7_2=K_3 & 0 & 7 & -2 & h_6  & h_3 \\ \hline
7_3 & 0 & 7 & -4 & h_4 & h_5  \\ \hline
7_4 & 0 & 7 & -2 & (h_4)^2 & h_3 \\ \hline
7_5 & 7 & 0 & 4 &  h_3 h_4  & h_4 \\ \hline
7_6 & 5 & 2 & 2 &  h_3 h_4 & h_3^2  \\ \hline
7_7 & 3 & 4 & 0 & h_3^3 & h_3^2 \\ \hline \hline
8_1=K_{-3} & 6 & 2 & 0 & h_3 & h_7 \\ \hline
8_2 & 6 & 2 & 4 & h_3 h_6 & h_3 \\ \hline
8_3 & 4 & 4 & 0 & h_5 & h_5 \\ \hline
8_4 & 4 & 4 & 2 & h_4 h_5  & h_3    \\ \hline
8_5 & 2 & 6 & -4 & h_3   &  ??? 
\\ \hline
\hline
K_p, p>0 & 0 & 2p+1 & -2 & h_{2p}^* & h_3  \\ \hline
K_p, p<0 & 2|p| & 2 & 0  & h_3 & h_{2|p|+1} \\ \hline
T(2,p),p>0 & 2p+1 & 0 & 2p & h_{2p+1} & 1 \\ \hline
\end{array}
$$
}
where, for a positive natural number $b$, $h_b$ are the unary 
theta and false theta series

\begin{equation*}
\lbl{eq.ha}
h_b(q)=\sum_{n \in \BZ} \ve_b(n) \, q^{\frac{b}{2} n(n+1)-n}
\end{equation*}
where
\begin{equation*}
\ve_b(n) =\begin{cases}
(-1)^n  & \text{if $b$ is odd} \\
1       & \text{if $b$ is even and $n \geq 0$} \\
-1       & \text{if $b$ is even and $n<0$} 
\end{cases}
\end{equation*}
Observe that 
$$
h_1(q)=0, \qquad h_2(q)=1, \qquad h_3(q)=(q)_\infty \,.
$$
In the above table, $c_+$ (resp. $c_-$) denotes the number of positive 
(resp., negative) crossings of an alternating knot $K$, and 
$\Phi^*_{K,0}(q)=\Phi_{-K,0}(q)$ denotes the $q$-series of the mirror $-K$ of $K$,
and $T(2,p)$ denotes the $(2,p)$ torus knot.

Concretely, the above table predicts the following identities 

{\small
\begin{eqnarray*}
(q)_\infty^{-2} &=& \sum_{a,b,c \geq 0}
(-1)^{a} \frac{q^{\frac{3}{2} a^2 + a b + a c + b c+
\frac{1}{2} a+ b+c}}{(q)_{a} (q)_{b} (q)_{c} (q)_{a+b}
(q)_{a+c}} \\ 
(q)_\infty^{-3} &=& \sum_{\substack{a,b,c,d,e  \geq 0 \\
a+b=d+e}} 
(-1)^{b+d}\frac{q^{\frac{b^2}{2}+ \frac{d^2}{2}+b c + a c + a
d+b e + \frac{a}{2}+c + \frac{e}{2} }}{(q)_{b+c}(q)_{a}
(q)_{b} (q)_{c} (q)_{d} (q)_{e} (q)_{c+d}} \\
(q)_\infty^{-4} &=& \sum_{\substack{a,b,c,d,e,f \geq 0 \\
a+e \geq b, b+f \geq a}}
(-1)^{a-b+e} 
\frac{q^{\frac{a}{2}+\frac{3 a^2}{2}+\frac{b}{2}+\frac{b^2}{2}+c+a c+d+a d+c d+\frac{e}{2}+2 a e-2 b e+d e+\frac{3 e^2}{2}-a f+b f+f^2}}{(q)_a (q)_b (q)_c (q)_{a+c} 
(q)_d (q)_{a+d} (q)_e (q)_{a-b+e} (q)_{a-b+d+e} (q)_f (q)_{-a+b+f}}
\end{eqnarray*}
}
corresponding to the knots
$$
3_1 \psdraw{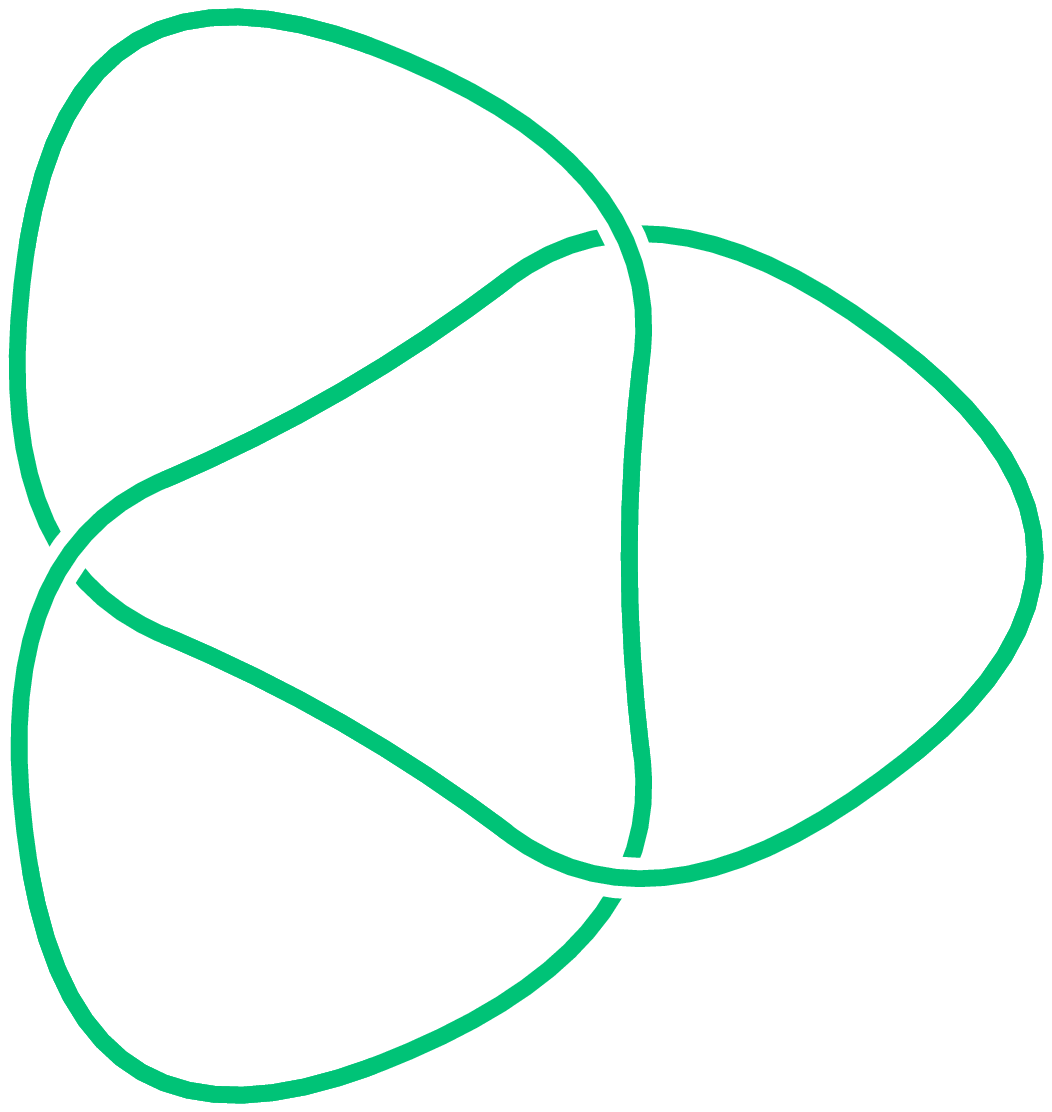}{1in}
\qquad\qquad 4_1 \psdraw{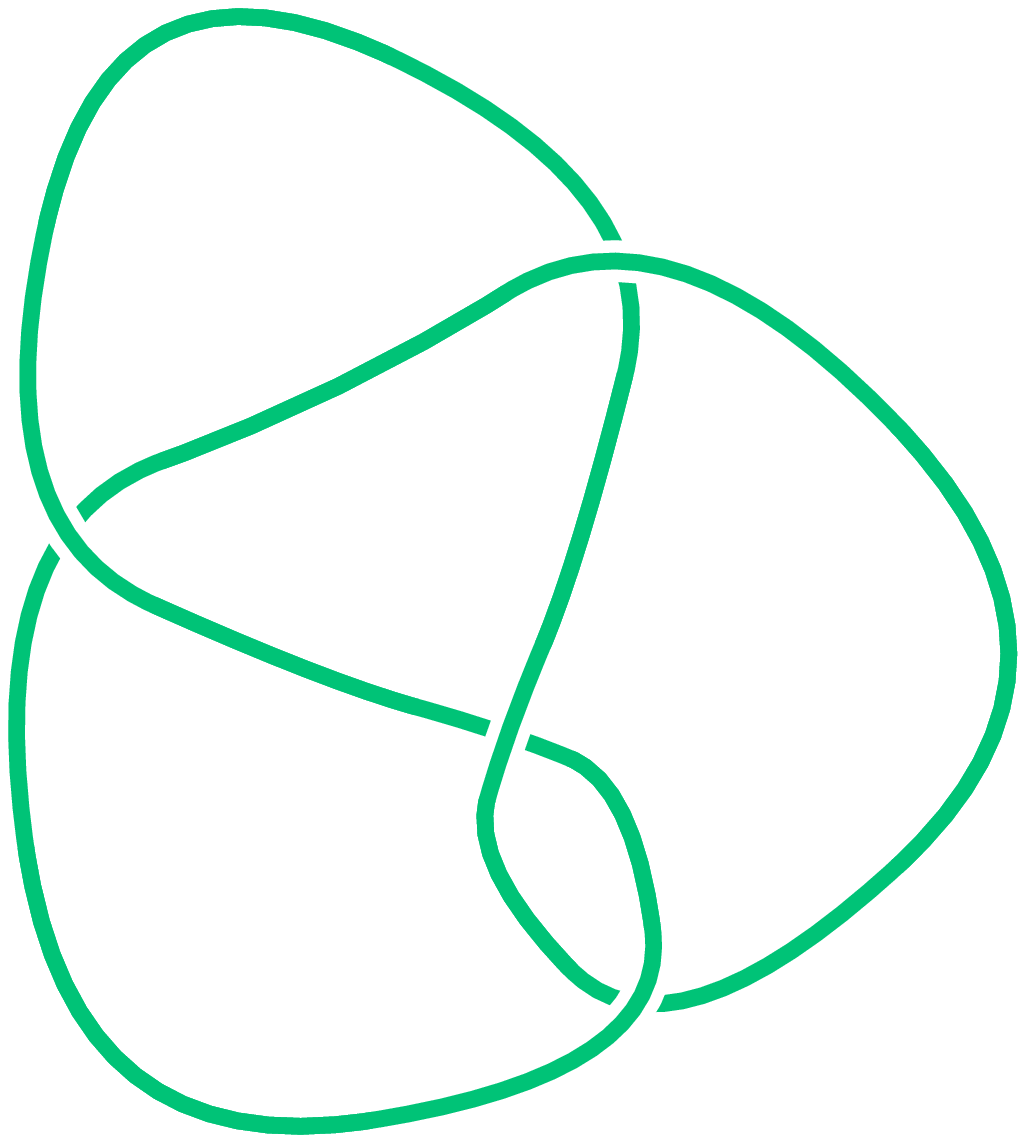}{1in}
\qquad\qquad 6_3 \psdraw{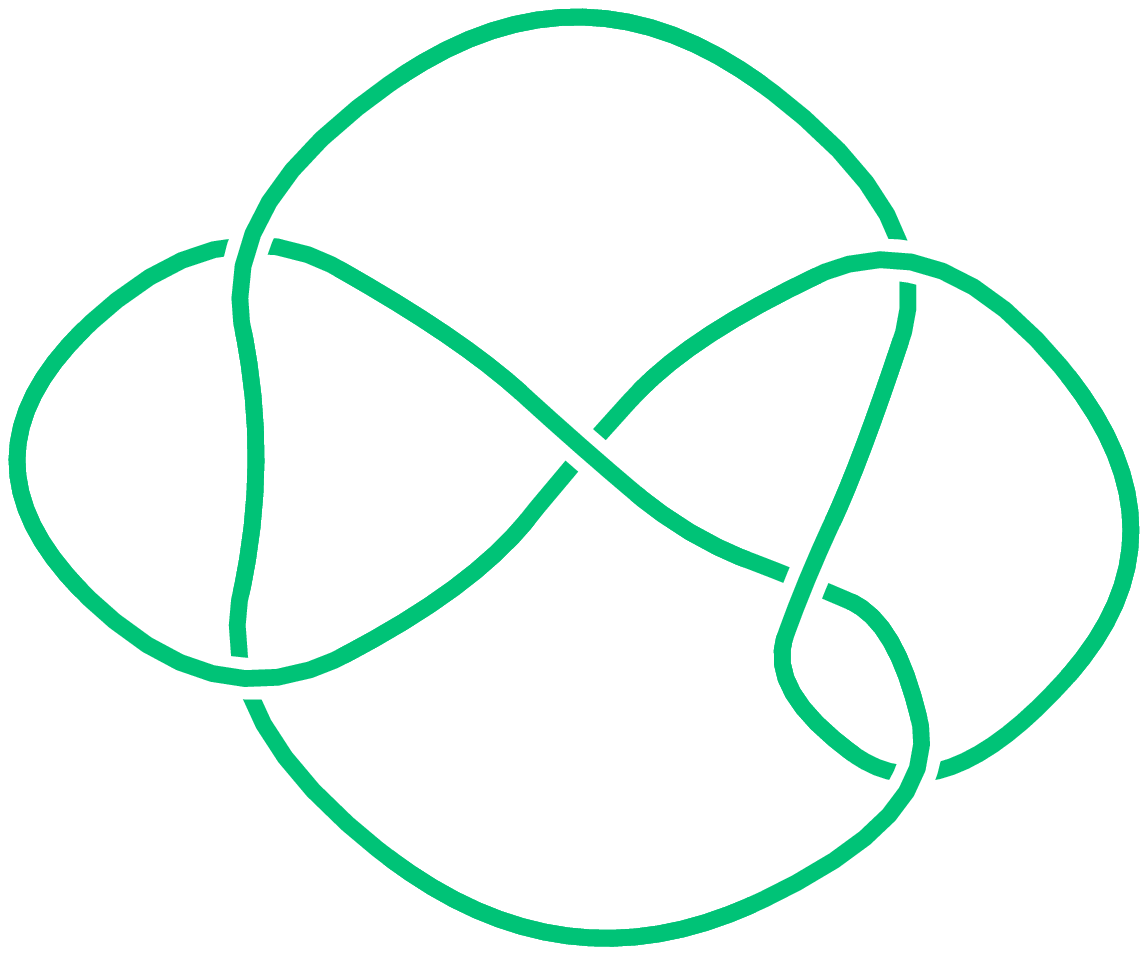}{1in}
$$

Some of the identities of the above table have been consequently 
proven \cite{ArDas}. In particular this settles the (mock)-modularity 
properties of the series $\Phi_{K,0}(q)$ for all but one knot.
The $q$-series of the remaining knot $8_5$ is given by an $8$-dimensional 
Nahn sum 

$$
\Phi_{8_5,0}(q)=
(q)_\infty^8 \sum_{\substack{a,b,c,d,e,f,g,h  \geq 0 \\ a+f \geq b}}
S(a,b,c,d,e,f,g,h)
\qquad\qquad
8_5 \psdraw{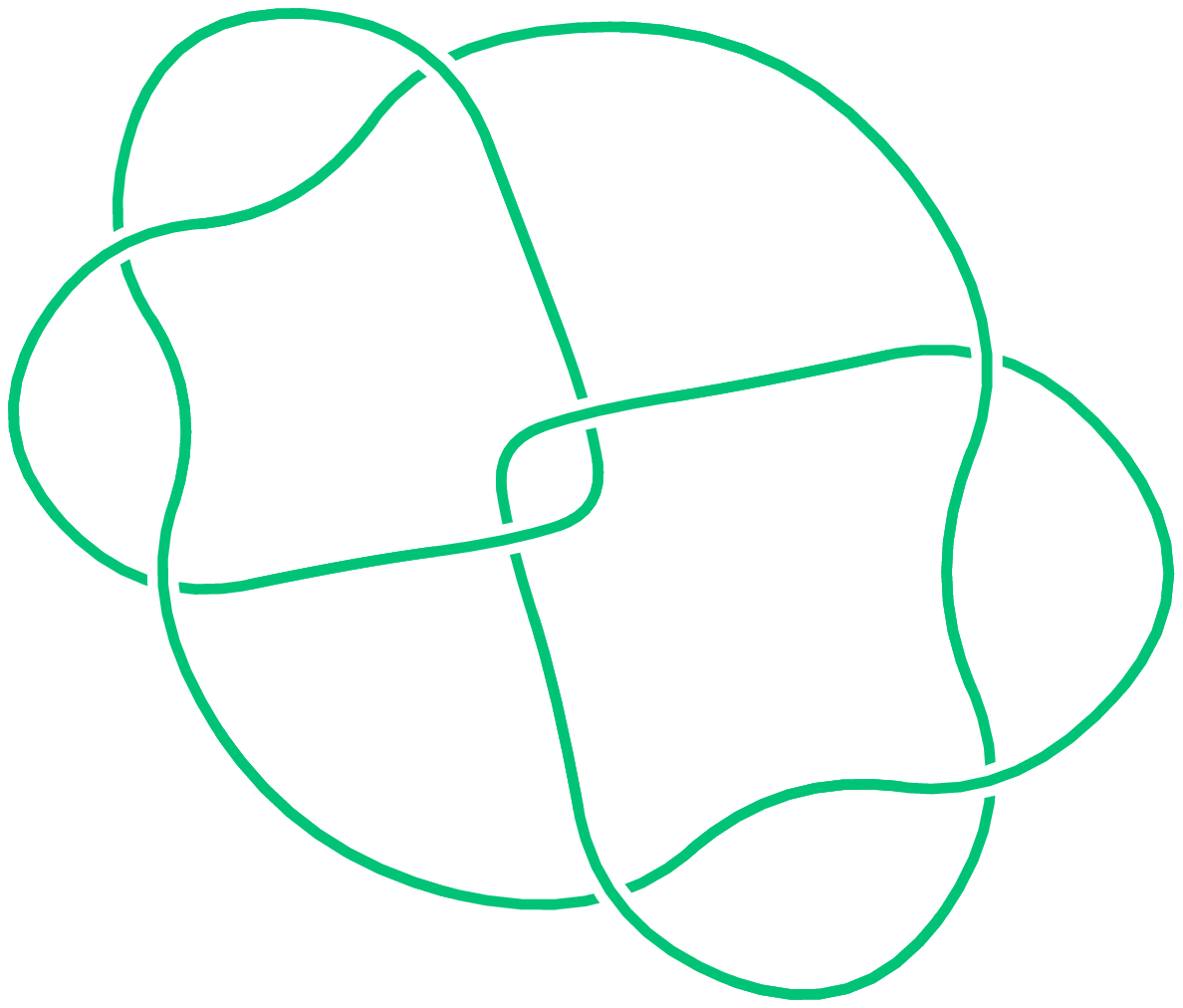}{1.2in}
$$
where $S=S(a,b,c,d,e,f,g,h)$ is given by

{\tiny
$$
S=
(-1)^{b+f} \frac{q^{2 a+3 a^2-\frac{b}{2}-2 a b+\frac{3 b^2}{2}+c+a c+d+a d+c d+e+a e+d e+\frac{3 f}{2}+4 a f-4 b f+e f+\frac{5 f^2}{2}+g+a g-b g+e g+f g+h+a h-b h+f h+g h}}{
(q)_{a} (q)_{b} (q)_{c} (q)_{d} (q)_{e} (q)_{f} (q)_{g} (q)_{h} 
(q)_{a+c} 
(q)_{a+d} 
(q)_{a+e} 
(q)_{a-b+f} 
(q)_{a-b+e+f} 
(q)_{a-b+f+g} 
(q)_{a-b+f+h}} \,.
$$
}
The first few terms of the series $\Phi_{8_5,0}(q)$, which somewhat simplify
when divided by $(q)_\infty$, are given by
\vspace{0.2cm}

\begin{math}
\Phi_{8_5,0}(q)/(q)_\infty=
1-q+q^2-q^4+q^5+q^6-q^8+2 q^{10}+q^{11}+q^{12}-q^{13}-2 q^{14}+2 q^{16}+3 q^{17}+2 q^{18}+q^{19}-3 q^{21}-2 q^{22}+q^{23}+4 q^{24}+4 q^{25}+5 q^{26}+3 q^{27}+q^{28}-2 q^{29}-3 q^{30}-3 q^{31}+5 q^{33}+8 q^{34}+8 q^{35}+8 q^{36}+6 q^{37}+3 q^{38}-2 q^{39}-5 q^{40}-6 q^{41}-q^{42}+2 q^{43}+9 q^{44}+13 q^{45}+17 q^{46}+16 q^{47}+14 q^{48}+9 q^{49}+4 q^{50}-3 q^{51}-8 q^{52}-8 q^{53}-5 q^{54}+3 q^{55}+14 q^{56}+21 q^{57}+27 q^{58}+32 q^{59}+33 q^{60}+28 q^{61}+21 q^{62}+11 q^{63}+q^{64}-9 q^{65}-11 q^{66}-11 q^{67}-2 q^{68}+9 q^{69}+27 q^{70}+40 q^{71}+56 q^{72}+60 q^{73}+65 q^{74}+62 q^{75}+54 q^{76}+39 q^{77}+23 q^{78}+4 q^{79}-9 q^{80}-16 q^{81}-14 q^{82}-3 q^{83}+16 q^{84}+40 q^{85}+67 q^{86}+92 q^{87}+114 q^{88}+129 q^{89}+135 q^{90}+127 q^{91}+115 q^{92}+92 q^{93}+66 q^{94}+35 q^{95}+9 q^{96}-12 q^{97}-14 q^{98}-11 q^{99}+13 q^{100}+O(q)^{101} \,.
\end{math}
\vspace{0.2cm}

We were unable to identify $\Phi_{8_5,0}(q)$ with a known $q$-series. Nor
were we able to decide whether it is a mock-modular form \cite{Za.mock}.
It seems to us that $8_5$ is not an exception, and that the mock-modularity
of the $q$-series $\Phi_{8_5,0}(q)$ is an open problem. 

\begin{question}
\lbl{que.mock}
Can one decide if a generalized Nahm sum is a mock-modular form?
\end{question}

\section{Modularity and Stability}
\lbl{sec.ms}

Modularity and Stability are two important properties of quantum knot
invariants. The Kashaev invariant $\la K \ra$ and the $q$-series $\Phi_{K,0}(q)$
of a knotted 3-dimensional object have some common features, namely 
asymptotic expansions at roots of unity approached radially 
(for $\Phi_{K,0}(q)$) and on the unit circle (for $\la K \ra$), 
depicted in the following figure
$$
\psdraw{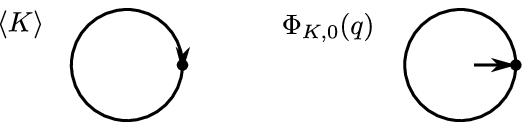}{2.5in}
$$
The leading asymptotic expansions of $\la K \ra$ and $\Phi_{K,0}(q)$ are 
governed by elements of the Bloch group as is the case of the Kashaev 
invariant and also the case of the radial limits of Nahm sums \cite{VZ}. 
In this section
we discuss a conjectural relation, discovered accidentally by Zagier and
the author in the spring of 2011, between the asymptotics of $\la 4_1\ra$ and
$\Phi_{6j,0}(q)$, where $6j$ is the $q$-6j symbol of the tetrahedron graph 
whose edges are colored with $2N$ \cite{Co,GV} 
$$
\psdraw{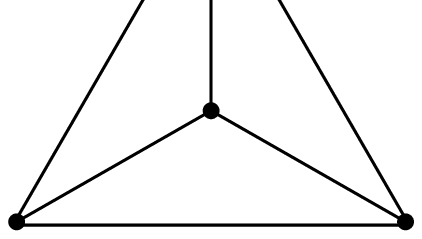}{0.5in}
$$
The evaluation of the above tetrahedron graph $J^+_{6j,N}(q) \in 1+q\BZ[q]$ 
is given explicitly by \cite{Co,GV}
$$
J^+_{6j,N}(q) =
\frac{1}{1-q} \sum_{n=0}^{N} (-1)^n
\frac{ q^{\frac{3}{2}n^2+\frac{1}{2}n}}{(q)_n^3}
\frac{(q)_{4N+1-n}}{(q)_n^3 (q)_{N-n}^4}\,.
$$
The sequence $(J^+_{6j,N}(q))$ is stable and the corresponding series
$F_{6j}(x,q)$ is given by
$$
F_{6j}(x,q)=
\frac{1}{(1-q)(q)_\infty^3}
\sum_{n=0}^\infty
(-1)^n
\frac{ q^{\frac{3}{2}n^2+\frac{1}{2}n}}{(q)_n^3}
\frac{(x q^{-n})_\infty^4}{(x^4 q^{-n+1})_\infty} \in \BZ((q))[[x]]\,,
$$
where as usual $(x)_\infty=\prod_{k=0}^\infty (1-xq^k)$ and 
$(q)_n=\prod_{k=1}^n (1-q^k)$. In particular,
$$
\lim_N J^+_{6j,N}(q)=\Phi_{6j,0}(q)= \frac{1}{(1-q)(q)_\infty^3}
\sum_{n=0}^\infty (-1)^n 
\frac{q^{\frac{3}{2} n^2 + \frac{1}{2}n}}{(q)_n^3}\,.
$$
Let 
$$
\phi_{6j,0}(q) = \frac{(q)^4_\infty}{1-q}\Phi_{6j,0}(q)= (q)_\infty
 \sum_{n=0}^\infty (-1)^n 
\frac{q^{\frac{3}{2} n^2 + \frac{1}{2}n}}{(q)_n^3} \,.
$$
The first few terms of $\phi_{6j,0}(q)$ are given by
\vspace{0.2cm}

\begin{math}
\phi_{6j,0}(q) =1-q-2 q^2-2 q^3-2 q^4+q^6+5 q^7+7 q^8+11 q^9+13 q^{10}+16 q^{11}+14 q^{12}+14 q^{13}+8 q^{14}-12 q^{16}-26 q^{17}-46 q^{18}-66 q^{19}-90 q^{20}-114 q^{21}-135 q^{22}-155 q^{23}-169 q^{24}-174 q^{25}-165 q^{26}-147 q^{27}-105 q^{28}-48 q^{29}+37 q^{30}+142 q^{31}+280 q^{32}+435 q^{33}+627 q^{34}+828 q^{35}+1060 q^{36}+O(q)^{37} \,.
\end{math}

\vspace{0.2cm}
The next conjecture which combines stability and modularity of two knotted
objects has been numerically checked around complex roots
of unity of order at most $3$.

\begin{conjecture}
\lbl{conj.416j}
As $X \longto +\infty$ with bounded denominator, we have
$$
\phi_{6j,0}(e^{-1/X}) = \phi_{4_1}(X)/X^{1/2} 
+ \overline{\phi_{4_1}(-\bar X)/(- \bar X)^{1/2}} \,.
$$
\end{conjecture}

\bibliographystyle{hamsalpha}\bibliography{biblio}
\end{document}